\documentclass[11pt]{amsart}

\usepackage{amsmath}
\usepackage{amsthm}
\usepackage{amssymb}
\usepackage{delarray}
\usepackage[dvips]{graphics}
\usepackage{epsfig}
\usepackage{color}

\renewcommand\bar\overline

\newtheorem{thm}{Theorem}[section]

\newtheorem{prop}[thm]{Proposition}
\newtheorem{cor}[thm]{Corollary}
\newtheorem{conj}[thm]{Conjecture}

\theoremstyle{definition}
\newtheorem{defn}[thm]{Definition}

\theoremstyle{remark}
\newtheorem{rmk}[thm]{Remark}
\newtheorem{obs}[thm]{Observation}

\title{Longest Induced Cycles in Cayley Graphs}
\author{Elena D. Fuchs \and Justin Sinz}
\address {University of California, Berkeley\\ Berkeley, CA \and University of Chicago\\ Chicago, IL}
\email{lenfuchs@berkeley.edu \and conrad@math.uchicago.edu}

\begin{document}

\maketitle

\begin{abstract}
In this paper we study the length of the longest induced cycle in the unitary Cayley graph $X_n = Cay(\mathbb Z_n; U_n)$, 
where $U_n$ is the group of units in $\mathbb Z_n$.  Using residues modulo the primes dividing $n$, we introduce a 
representation of the vertices that reduces the problem to a purely combinatorial question of comparing strings of 
symbols.  This representation allows us to prove that the multiplicity of each prime dividing $n$, and even the 
value of each prime (if sufficiently large) has no effect on the length of the longest induced cycle in $X_n$.  
We also see that if $n$ has $r$ distinct prime divisors, $X_n$ always contains an induced cycle of length $2^r+2$, 
improving the $r \ln r$ bound of Berrezbeitia and Giudici.  Moreover, we extend our results for $X_n$ to conjunctions 
of complete $k_i$-partite graphs, where $k_i$ need not be finite, and also to unitary Cayley graphs on any quotient of a 
Dedekind domain.

\end{abstract}

\section{Introduction}\label{intro}

For a positive integer $n$, let the unitary Cayley graph $X_n = Cay(\mathbb Z_n , U_n)$ be defined as follows:

(1) The vertex set of $X_n$, denoted by $V(n)$, is $\mathbb Z_n$, the ring of integers modulo $n$.

(2) The edge set of $X_n$ is denoted by $E(n)$, and, for $x, y \in V(n)$, $\{x,y\} \in E(n)$ if and only 
if $x-y \in U_n$, where $U_n$ is the set of units in the ring $\mathbb Z_n$.

The central problem adressed in this paper is to find the length of the longest induced cycle in $X_n$.  This problem 
was first considered by Berrizbeitia and Giudici \cite{OrigArt}, who were motivated by its applications to 
chromatic uniqueness.

Throughout the paper, we let $n= p_1 ^{a_1}p_2^{a_2} \dots p_r^{a_r}$, where the $p_i$ are distinct primes, and 
$a_i \geq 1$.  Then we denote the length of the longest induced cycle in $X_n$ by $M(n)$.  We let $m(r)= \max_{n}M(n)$, 
where the maximum is taken over all $n$ with $r$ distinct prime divisors.  In \cite{OrigArt}, Berrizbeitia and 
Giudici bound $m(r)$ by $$r \ln r \leq m(r) \leq 9r!.$$  A simple change to the proof of the upper bound 
provided in \cite{OrigArt} yields the better upper bound of $m(r) \leq 6r!$.

Our goal is to determine better bounds for $m(r)$, as well as to extend what we find to other graphs.  
In Section~\ref{ResRep}, we introduce a useful representation of the vertices in $X_n$ according to their residues modulo 
the prime divisors of $n$.  This representation immediately yields several helpful properties of the longest induced 
cycles in these graphs.  In particular, we prove that we can disregard the multiplicities of the prime divisors of 
$n$, so we can reduce our problem to square-free $n$.  Also, we show 
that $M(n)$ depends only on $r$, and in fact $M(n) = m(r)$ as long as the primes dividing $n$ are all large enough.  
In Section~\ref{Bound}, we use the vertex representation introduced in Section~\ref{ResRep} to construct an induced 
cycle of length $2^r +2$ in the graph $X_n$, where $n$ has $r$ distinct prime divisors, thus raising the lower 
bound on $m(r)$ substantially.  We also note that this construction is valid for any $n$, no matter what its prime 
divisors are, so this provides a lower bound for $M(n)$.  Section~\ref{general} 
contains a generalization of our results to conjunctions of complete $k_i$-partite graphs, as well as to unitary 
Cayley graphs on products of local rings, which include the unitary Cayley graphs on Dedekind rings.  We conclude with 
open questions that we believe may be solved with the use of the vertex representation that we introduce in 
Section~\ref{ResRep}.

\section{Residue Representation}\label{ResRep}

Recall that $n = p_1^{a_1} p_2^{a_2} \cdots p_r^{a_r}$, where the $p_i$ are prime.  We will represent the vertices 
of $X_n$ in a way that will reduce the process of finding induced cycles in $X_n$ to checking for similarities between 
strings of numbers in an array.

It is clear that the following is equivalent to the definition of $E(n)$ in the introduction:

\begin{obs}\label{edge}
For $x,y \in V(n)$, we have that $ \{ x,y \} \in E(n)$ if and only if 
$$x \not\equiv y \pmod {p_i} \mbox{, for all } 1 \leq i \leq r.$$
Likewise, $ \{ x,y \} \not\in E(n)$ if and only if
$$x \equiv y \pmod {p_i} \mbox{, for some } 1 \leq i \leq r.$$
\end{obs}

So, in fact, to know whether $x$ and $y$ are adjacent we need only their residues modulo the primes $p_i$.  With this 
in mind, we introduce the following representation of the vertices:

\begin{defn}\label{vrr}
\item{(i)} Let $x \in V(n)$, such that 
$$x \equiv \alpha_i \pmod {p_i} \mbox {, where } 1 \leq i \leq r \mbox { and } 0 \leq \alpha_i < p_i.$$
We then define the \emph {residue representation} of $x$ to be the unique string $\alpha_1 \alpha_2 \cdots \alpha_r$, where 
$\alpha_k$ is the $k$\emph{th term}, and we 
write $x \approx \alpha_1 \alpha_2 \cdots \alpha_r$.
\item{(ii)} Let $x , y \in V(n)$.  If the $k$th term of the residue representation of $x$ is the same as the $k$th term 
of the residue representation of $y$, we say that $x$ has a \emph{similarity} with $y$.
\end{defn}

Combining Observation~\ref{edge} and Definition~\ref{vrr}, vertices $x, y \in V(n)$ are adjacent if and only if $x$ 
has no similarities with $y$.  So, in fact, the only property of the residues modulo $p_i$ that we use in constructing 
induced cycles is that they form a set of size $p_i$, and we verify that a subgraph is an induced cycle by checking 
that consecutive vertices do not have any similarities, and that any pair of non-consecutive vertices has at least one 
similarity.

Also, we note that for $n$ not square-free, a string may be the residue representation of multiple vertices.  For example, 
if $n = 12$, both $0$ and $6$ have residue representation $00$.  However, the adjacency of vertices depends 
only on their residue representations, and, by the Chinese Remainder Theorem, every string represents at least one vertex.

This representation greatly simplifies inspection of induced cycles.  In fact, we can extend residue representation for a 
vertex to any induced subgraph:

\begin{defn}

\item {(i)} Let $S$ be an induced subgraph of $X_n$, where $V(S) = (v_0, v_2, \dots ,v_{k-1})$, with 
$v_i \approx \alpha_{i1} \alpha_{i2} \cdots \alpha_{ir}$, and $0 \leq i \leq k-1$.  
We then define the \emph {residue representation of} $S$ to be the array 

\[
\begin{array}{cccc}
\alpha_{01} & \alpha_{02} & \cdots & \alpha_{0r}\\
\alpha_{11} & \alpha_{12} & \cdots & \alpha_{1r}\\
\vdots & \vdots & & \vdots\\
\alpha_{(k-1)1} & \alpha_{(k-1)2} & \cdots & \alpha_{(k-1)r}.\\
\end{array}
\]

\item{(ii)} The \emph {residue set} of $S$ is the set of residues 
\begin{equation*}
\bigcup_{\substack{0 \leq i \leq k-1\\ 1 \leq j \leq r}}\ \{\alpha_{ij} \}
\end{equation*}
used in its residue representation.

\end{defn}

So, if an induced subgraph $S$ is a $k$-cycle in $X_n$, we can permute the rows of the residue representation of $S$
so that the $i$th row has a similarity with the $j$th row if and only if $i-j \not\equiv \pm 1 \pmod k$.  
Figure~\ref{eg} displays the residue representation of an induced $6$-cycle for $r=2$ and for $r=3$.

\begin{figure}
\[
\begin{array}{cccccc}
0 & 0 &\hspace{1 in} &0&0&0\\
1 & 1 &\hspace{1 in} &1&1&1\\
0 & 2 &\hspace{1 in} &0&0&2\\
1 & 0 &\hspace{1 in} &1&2&0\\
0 & 1 &\hspace{1 in} &0&0&1\\
1 & 2 &\hspace{1 in} &1&1&2\\
\end{array}
\]
\caption{In these residue representations of an induced $6$-cycle for $r=2$ on the left, and for $r=3$ on the right, 
it is easy to see that two consecutive rows (including the $1$st and $6$th rows) have no similarities, and any two 
non-consecutive rows have at least one similarity.  The residue set for each cycle is $\{ 0, 1, 2 \}$.}
\label{eg}
\end{figure}

An important property of an induced cycle of length greater than $4$ is that it cannot contain two vertices with the same 
residue representation.

\begin{prop}\label{nosame}
The residue representation of a $k$-cycle $C$, with $k>4$, cannot contain two identical rows.
\end{prop}

\begin{proof}
Suppose there are two vertices $x$ and $y$ in $C$ that have the same residue representation.  Then a vertex $z$ of $C$ has 
no similarity with $x$ if and only if it has no similarity with $y$, meaning that $x$ and $y$ have precisely the same 
neighbors in $C$.  However, a vertex in an induced cycle is adjacent to exactly two other vertices in the cycle, so 
$C$ can have at most $4$ vertices, contradicting $k>4$.  Thus the residue representation of $C$ cannot contain two 
identical rows.
\end{proof}

It is important that, once we have written an induced cycle in terms of its residue representation, we can permute 
the residues in each column to obtain an induced cycle of equal length.

\begin{obs}\label{permute}
Let the $j$th column in the residue representation of an induced $k$-cycle $C$ in $X_n$ be 
\[
\begin{array}{c}
\alpha_{0j}\\
\alpha_{1j}\\
\vdots\\
\alpha_{(k-1)j},\\
\end{array}
\]
and suppose this column contains $l_j$ distinct residues, $\{a_1, a_2, \dots, a_{l_j}\}$.  Then let $\pi$ be a 
permutation of $\{a_1, a_2, \dots, a_{l_j}\}$, and replace the $j$th column of $C$ by 
\[
\begin{array}{c}
\pi(\alpha_{0j})\\
\pi(\alpha_{1j})\\
\vdots\\
\pi(\alpha_{(k-1)j}).\\
\end{array}
\]
We then have a new induced $k$-cycle in $X_n$, since we have not changed the similarities between any of the rows in $C$.
\end{obs}

We now use the Observation~\ref{permute} to define isomorphisms between induced $k$-cycles in $X_n$.

\begin{defn}\label{iso}
Two induced $k$-cycles, $C$ and $C'$, are called \emph{isomorphic} if the $j$th column of the residue representation of 
$C'$ is obtained by permuting the residues in the $j$th column of $C$, as described in Observation~\ref{permute}.
\end{defn}

Note that the first two rows in Figure~\ref{eg} are $000$ and $111$.  Because of this, all of the rows that 
are not adjacent to either of the first two have to contain both a $0$ and a $1$.  Similarly, the third row 
in the cycle must contain a $0$, and the last row in the cycle must contain a $1$.  This is a useful criterion 
for induced cycles in general.

\begin{rmk}\label{firsttwo}
Any induced cycle $C$ in $X_n$ is isomorphic to an induced cycle $C'$ of the same length so that the first two rows in 
the residue representation of $C'$ are $00 \cdots 0$ and $11 \cdots 1$.
\end{rmk}

In order to obtain such a $C'$, we need only to map the first two elements in every column of $C$ to $0$ and $1$, 
respectively.  Note that the first two elements in each column are always different -- if they were not, the first 
and the second row in the residue representation of $C$ would have a similarity, which contradicts their adjacency.

This tells us that all but four of the rows in our induced cycles will have to contain both a $0$ and a $1$, which may 
limit the residue sets and consequently the lengths of the cycles.

Another interesting fact that becomes evident with the use of residue representation is the following proposition.

\begin{prop}\label{incr}
The value $m(r)$ increases with $r$.  Specifically, if $X_n$ contains an induced cycle of length $k$, and $q > 2$ 
is a prime not dividing $n$, then $X_{qn}$ also contains a cycle of length $k$.  If $k$ is even, we can also 
allow $q = 2$.
\end{prop}

\begin{proof}
Let $n= p_1^{a_1} p_2^{a_2} \cdots p_r^{a_r}$, where the exponents $a_i$ are positive integers, and $p_i$ are 
distinct primes.  Suppose $X_n$ contains an induced cycle $C$ of length $k$.  We denote the residue representations 
of the vertices of $C$ by $v_0, v_1, \dots , v_{k-1}$, where each $v_i$ is a string of length $r$.  Let 
$n' = q n$, where $q \not = 2$ is prime, $q \not = p_i \mbox { for all } 1 \leq i \leq r$.  Then we will 
show that $X_{n'}$ also contains a cycle of length $k$ by constructing an induced cycle $C'$ in 
$X_{n'}$, denoting the residue representations of the vertices of $C'$ by $w_0, w_1, \dots , w_{k - 1}$.

If $k$ is even, let $w_i = 0v_i$ for even $i$, and let $w_i = 1v_i$ for odd $i$.  Notice that we do not introduce any 
similarities between two rows that were adjacent in $C$, so two consecutive rows in $C'$ are adjacent, as desired.  
Similarly, if $\{v_i, v_j\} \not\in E(n)$, they have a similarity, say, in the $l$th term.  Then $w_i$ and $w_j$ have 
a similarity in the $(l+1)$st term, and so $\{w_i, w_j\} \not\in E(n')$.  Thus we introduce no new adjacencies in the 
construction of $C'$, so $C'$ is indeed an induced $k$-cycle in $X_{n'}$.

If $k$ is odd, let $w_i = 1v_i$ for odd $i$, let $w_i  = 0v_i$ for even $i \not = k - 1$, and let 
$w_{k - 1} = 2v_{k - 1}$ (this is possible since $q \not = 2$).  Again, we note that we do not introduce any 
similarities between two rows that were adjacent in $C$, so two consecutive rows in $C'$ are adjacent, as desired.  
Also, if $\{ v_i, v_j \} \not\in E(n)$, we have that $\{ w_i, w_j \} \not\in E(n')$ by the argument above.  Thus we 
introduce no new adjacencies in the construction of $C'$, so $C'$ is indeed an induced $k$-cycle in $X_{n'}$.

By starting with a cycle $C$ in $X(n)$ that has length $m(r)$, we see that $m(r+1) \geq m(r)$, as desired.  
\end{proof}

\begin{cor}\label{greatsix}
If $r \geq 2$, and $n$ is square-free, then $M(n) \geq 6$.
\end{cor}

\begin{proof}
For $r = 2$, we have constructed a $2$-cycle of length $6$ in Figure~\ref{eg}, so $m(2) \geq 6$.  
Proposition~\ref{incr} shows that $m(r)$ is nondecreasing, so we have that, if $r > 2$, $m(r) \geq m(2) \geq 6$, as 
desired.
\end{proof}

We now prove that, in calculating $M(n)$, we need consider only those $n$ that are square-free.

\begin{thm}\label{Mult}
For $n = p_1^{a_1}p_2^{a_2}\cdots p_r^{a_r}$, 
and $n' = p_1p_2\cdots p_r$, where $r \not= 1$, $M(n) = M(n')$.
\end{thm}

\begin{proof}
(1) First we show that $M(n) \geq M(n')$.  In particular, we show $X_n$ contains cycles of length $M(n')$.  
Note that since $n$ and $n'$ have the same prime divisors, if $x, y < n$, then $x-y \in U_n$ if and only if 
$x-y \in U_{n'}$.  So, in particular, the induced subgraph of $X_n$ on vertices $0, 1, \dots , n'-1$ is 
precisely $X_{n'}$.  Thus any induced cycle on $X_{n'}$ can be mapped to an induced cycle in 
$\{0, 1, \dots , n'-1 \} \subset X_n$, and so there is an induced cycle of length $M(n')$ in $X_n$, as desired.

(2) Now we show that $M(n) \leq M(n')$, or that there is no induced cycle of length greater than $M(n')$ in $X_n$.  
Since $n'$ is square-free, Corollary~\ref{greatsix} implies that $M(n') \geq 6$.  Suppose there is an induced cycle, 
$C_l$, of length $l > M(n')$ in $X_n$.  Then, in particular, $l > 6$.  Using residue representation, write $C_l$ in 
terms of residues $\pmod {p_1, p_2, ... , p_r}$.  If no two vertices 
in $C_l$ are denoted by the same string of residues, then we can view the residue representation of $C_l$ as a residue 
representation of an induced $l$-cycle in $X_{n'}$.  Since $l > M(n')$, this contradicts the assumption that $M(n')$ 
is the maximum length of an induced cycle in $X_{n'}$.  Thus there exist two vertices in $C_l$ 
that have identical residue representations.  However, by Proposition~\ref{nosame}, this means $l \leq 4$, 
contradicting the previous deduction that that $l > 6$.  We conclude that, indeed, there are no induced cycles of length 
$l > M(n')$ in $X_n$.
\end{proof}

\begin{prop}\label{Multprop}
Let $n' = p$, and $n = p^a$ where $p$ is a prime and $a > 1$.  Then $M(n') = 3$, and $M(n) = 4$.  So, $m(1) = 4$.
\end{prop}
\begin{proof}
Since the only non-unit in $\mathbb Z_p$ is $0$, $X_{n'}$ is a complete graph on $p$ vertices, and the 
longest induced cycle in $X_{n'}$ must hence have length $3$.  From Part (2) of the proof of 
Theorem~\ref{Mult}, we deduce that $M(n) \leq 4$.  In fact, $M(n) = 4$, since the subgraph 
$(0, 1, p, p+1)$ is an induced cycle in $X_n$.
\end{proof}

\begin{prop}\label{mnmr}
For $n = p_1 ^{a_1} p_2 ^ {a_2} \cdots p_r ^{a_r}$ where the $p_i$ are large, $M(n) = m(r)$.
\end{prop}

\begin{proof} 
Since $M(n)$ depends only on the residues available to use in a residue representation of cycles.  If $n$ and $n'$ 
each have $r$ distinct prime divisors, $M(n') = m(r)$, and the residue representation of some longest induced cycle 
in $X_{n'}$ is also the residue representation of a subgraph $S$ of $X_n$, then $S$ will in fact be an induced cycle 
in $X_n$, so $M(n') = m(r)$.  Thus, as long as the prime divisors of $n$ yield enough residues for a residue 
representation of the longest cycle in $X_{n'}$, where $M(n') = m(r)$, we will have $M(n) = m(r)$.
\end{proof}

\section{A Lower Bound on $m(r)$}\label{Bound}

One important asset of introducing residue representation is that it gives us a way to construct a good lower bound 
on $m(r)$; we achieve the following lower bound as our main result in this section.

\begin{thm}\label{Lowbd}
For all positive integers $n$ with $r > 1$ distinct prime divisors, we have $m(r) \geq 2^r + 2$.
\end{thm}

In this section, we construct an induced subgraph of $X_n$ with $2^r + 2$ vertices, where $r$ is the number of 
distinct prime divisors of $n$, and provide two specific cycles produced by this construction.  We will then prove that 
this subgraph is indeed a cycle, and thus show that Theorem~\ref{Lowbd} holds.

In order to construct an  induced $2^r+2$-cycle in $X_n$, where $n = p_1p_2 \cdots p_r$, we first introduce some definitions, which are 
discussed in detail in \cite{CombAlg}, p. 433.

(i) An \emph{$n$-bit Gray Code} is an ordered, cyclic sequence of the $2^n$ $n$-bit binary strings called 
\emph{codewords}, such that successive codewords differ by the 
complementation of a single bit, and the starting codeword is taken to be $(00\cdots0)$.  We write this sequence 
in the form of a matrix, as shown below.

(ii) A \emph{Reflective Gray Code} (RGC) is defined recursively as follows: A $1$-bit RGC is merely the $2 \times 1$ matrix 
$0 \choose 1$.  If an $r$-bit RGC is the $2^r \times r$ binary matrix
\[
\left(
\begin{array}{l}
G_0\\
G_1\\
\vdots\\
G_{2^r - 1},
\end{array}
\right)
\]
then we define the $(r+1)$-bit RGC to be the $2^{r+1} \times {(r+1)}$ binary matrix
\[
\left(
\begin{array}{l}
0G_0\\
0G_1\\
0G_2\\
\vdots\\
0G_{2^r - 1}\\
1G_{2^r - 1}\\
1G_{2^r - 2}\\
\vdots\\
1G_1\\
1G_0
\end{array}
\right).
\]

Henceforth, we fix $r$ and index the codewords by $0, 1, \dots, 2^r -1$ $\pmod {2^r}$, denoting the $i$th codeword in an $r$-bit RGC 
by $G_i$, and the $i$th codeword in a $k$-bit RGC, where $k \not = r$, by $G_i^{(k)}$.  

(iii) The \emph{flip bit} in the $j$th codeword of a RGC is the position of the one bit that has changed from the 
$(j-1)$st codeword.

We will construct an induced subgraph of $X_n$ whose residue representation consists of the rows 
$v_0, v_1, \dots , v_M$, where $M = 2^r +1$, and $\{v_i , v_j\} \in E$ if and only if $i-j \equiv \pm 1 \pmod {2^r + 2}$.  
Let $v_{M-1} \approx 0100 \cdots 0$, and $v_M \approx 122 \cdots 2$.  
We define the rows $\{v_i : i \mbox{ even}, i \not= M-1\}$ by using the first half of an 
$r$-bit RGC with a slight modification.  Let $\widehat G_i$, for $i \not = 0$ be the $i$th codeword $G_i$ in an $r$-bit 
RGC, with the flip bit replaced by a $2$.  Let $\widehat G_0 = G_0$.  Then we define the even-indexed rows as follows: 
$v_{2i} = \widehat G_i$, for $0 \leq i < 2^{r-1}$.

We define the odd-indexed rows as follows: for $0 \leq i \leq 2^{r-1}$, let $v_{2i+1} = \bar{G_i}$, the complement 
of $G_i$.  So the subgraph we have constructed is 
$\{ \widehat G_0, \bar G_0, \widehat G_1, \dots , \widehat G_{2^{r-1} - 1},$ $\bar G_{2^{r-1} - 1}, v_{M-1}, v_M \}$.  
This gives us a subgraph consisting of $(2^r +2)$ vertices.

In Figure~\ref{cycles}, we display this construction for $r = 3$ and $r = 4$.

\begin{figure}
\[
\begin{array}{cccccccc}
0&0&0&\hspace{1 in} & 0&0&0&0\\
1&1&1&\hspace{1 in} & 1&1&1&1\\
0&0&2&\hspace{1 in} & 0&0&0&2\\
1&1&0&\hspace{1 in} & 1&1&1&0\\
0&2&1&\hspace{1 in} & 0&0&2&1\\
1&0&0&\hspace{1 in} & 1&1&0&0\\
0&1&2&\hspace{1 in} & 0&0&1&2\\
1&0&1&\hspace{1 in} & 1&1&0&1\\
0&1&0&\hspace{1 in} & 0&2&1&0\\
1&2&2&\hspace{1 in} & 1&0&0&1\\
 & & &\hspace{1 in} & 0&1&1&2\\
 & & &\hspace{1 in} & 1&0&0&0\\
 & & &\hspace{1 in} & 0&1&2&1\\
 & & &\hspace{1 in} & 1&0&1&0\\
 & & &\hspace{1 in} & 0&1&0&2\\
 & & &\hspace{1 in} & 1&0&1&1\\
 & & &\hspace{1 in} & 0&1&0&0\\
 & & &\hspace{1 in} & 1&2&2&2\\
\end{array}
\]
\caption{We construct two cycles using residue representation and our lower bound construction.  
On the left is an induced $10$-cycle for the graph 
$X_n$, where $n$ has three prime divisors ($r = 3$).  On the right is an induced $18$-cycle for the graph $X_n$, where $n$ 
has four prime divisors ($r = 4$).  Note that the rows in both cycles are derived as described from a $3$-bit Reflective 
Gray Code and a $4$-bit Reflective Gray Code, respectively.}
\label{cycles}
\end{figure}

To prove Theorem~\ref{Lowbd}, we must show that the subgraph we have constructed is indeed an induced cycle.  This can be 
reduced to showing that the following properties hold.

(i) Vertex $v_k$ is adjacent to $v_l$ if $k-l \equiv \pm 1 \pmod {2^r + 2}$.  In other words, 
$\{v_0, v_1, \dots v_M\}$ is a cycle.

(ii) If neither $k$ nor $l$ equals $M-1$ or $M$, and $|k-l| > 1$, then $v_k$ is not adjacent to $v_l$.

(iii) Vertex $v_M$ is not adjacent to $v_l \mbox { for } i \not= 0, M-1$, and vertex $v_{M-1}$ is not adjacent to 
$v_l \mbox { for } i \not= M-2, M$.

\begin{proof}
[Proof of Theorem~\ref{Lowbd}]

\item{(i)} First we show that any two consecutive rows among ${v_0, v_1, \dots , v_{M-2}}$ correspond to adjacent vertices. 
Among these rows, no odd-indexed row contains a $2$, and an even-indexed row $v_{2i}$ is merely the complement of $v_{2i+1}$ 
with one bit replaced by a $2$.  Thus every odd-indexed row among ${v_0, v_1, \dots , v_{M-2}}$ has no similarities 
with the row immediately above it.  Also, since any two consecutive codewords $G_i$ and $G_{i+1}$ 
in an $r$-bit RGC differ only in the flip bit of $G_{i+1}$, the codeword $\bar G_i$ differs from $G_{i+1}$ 
everywhere except in the flip bit.  However, in modifying $G_i$ to $\widehat G_i$ for $0 \leq i < 2^{r-1}$, we have replaced every flip bit by a 
$2$, so $v_{2i+1} = \bar G_i$, (which will contain no $2$'s), will differ completely from 
$v_{2i+2} = \widehat G_{i+1}$ if $i \not = 2^{r-1}-1$.  Thus every odd-indexed row among ${v_0, v_1, \dots , v_{M-4}}$ 
is adjacent to the row immediately below it.

It remains to show that $v_M$ is adjacent to $v_{M-1}$, that $v_M$ is adjacent to $v_0$ (these two claims are trivial 
by inspection), and that $v_{M-2}$ is adjacent to $v_{M-1}$.  Note that $v_{M-1}$ is precisely $G_{2^{r - 1}-1}$, 
since, by definition, 
$$G_{2^{r - 1}-1} = 0G_{2^{r - 2} - 1} ^ {(r-1)} = 01G_0^{(r-2)} = 0100 \cdots 0.$$  
Also, $v_{M-2}$ is, by definition, $\bar G_{2^{r - 1}-1}$.  Thus, indeed, $v_{M-2}$ is adjacent to $v_{M-1}$, and we 
have that $\{v_0, v_1, \dots v_M\}$ is a cycle.

\item{(ii)} It is trivial to show that no two rows whose indices have the same parity are adjacent, since all 
even-indexed rows begin with a $0$ and are thus not adjacent to each other, while all odd-indexed rows begin 
with a $1$ and are also not adjacent to each other.

Now, take an even-indexed row $v_{2i}$, with $0 \leq i < 2^{r-1}$, and an odd-indexed row $v_{2j+1}$, with 
$0 \leq j < 2^{r-1}$, such that $i \not= j$ and $i \not= j + 1$.  Suppose for the sake of contradiction that $v_{2i}$ 
is adjacent to $v_{2j+1}$.

By definition, $v_{2j+1} = \bar G_j$, $v_{2i} = \widehat G_i$, and $i \not= j$ by assumption.  By the definition 
of a RGC, $G_j$ and $G_i$ differ in at least one bit.  Since $i-j \not\equiv 1 \pmod {2^r}$, then $G_j$ and $G_i$ 
must differ in a bit that is not a flip bit for $G_i$.  Therefore $v_{2j+1} = \bar G_j$ will have at least one 
similarity with $v_{2i} = \widehat G_i$, and so $v_{2i}$ and $v_{2j+1}$ are not adjacent, contrary to our supposition.

So, indeed, if neither $k$ nor $l$ equals $M-1$ or $M$, and $|k-l| > 1$, then $v_k$ is not adjacent to $v_l$.

\item{(iii)} Since $v_M$ begins with a $1$, it is not adjacent to any of the odd-indexed rows, which also all begin with a $1$.  
Similarly, because all of the even-indexed rows except $v_0$ and $v_{M-1}$ have a $2$ in some spot after the initial $0$, 
and will thus have a similarity with $v_M \approx 122 \cdots 2$, no even-indexed row except $v_0$ and $v_{M-1}$ 
will be adjacent to $v_M$.

Since $v_{M-1}$ begins with a $0$, it is not adjacent to any of the even-indexed rows, which all begin with a $0$ as well.  
Also, note that $v_{M-2} = v_{2^r -1} = \bar G_{2^{r-1}-1} = 1011 \cdots 1$ is the complement of $v_{M-1}$, and that all 
odd-indexed rows except $v_M$ are distinct and contain only $0$'s and $1$'s.  Thus all odd-indexed rows except $v_M$ 
either complement or have a 
similarity with $V_{M-1} = 0100 \cdots 0$.  So all odd-indexed rows except for $v_{M-2}$ and $v_M$ are not adjacent to 
$v_{M-1}$.

Thus we have that vertex $v_M$ is not adjacent to $v_i \mbox { for } i \not= 0, M-1$, and 
vertex $v_{M-1}$ is not adjacent to $v_i \mbox { for } i \not= M-2, M$.
\end{proof}

Note that, for any $n = p_1p_2 \cdots p_r$, where $p_1 < p_2 < \cdots < p_r$ are primes, the cycle constructed above 
does not depend on the choice of $p_i$.  The first column of the cycle's residue representation contains residues 
$0$ and $1$ only, allowing for $p_1 =2$, and the residue set of the cycle is $\{0, 1, 2\}$, which puts no bounds 
on the rest of the primes $p_i$.

Also, Theorem~\ref{Mult} implies that our construction of a $(2^r +2 )$-cycle for $n' = p_1p_2...p_r, r > 1$ holds for 
$n = {p_1}^{a_1}{p_2}^{a_2}...{p_r}^{a_r}$, while Proposition~\ref{Multprop} implies that the lower bound in 
Theorem~\ref{Lowbd} holds for $r=1$.

\section{Generalizing to Other Graphs}\label{general}

A natural question to ask is what properties of the Cayley graph $X_n$ are necessary to obtain the results we have.  
It is noted in \cite{OrigArt} that, for $p$ prime and $a$ a positive integer, $X_{p^a}$ is complete $p$-partite.  
In fact, this tells us that for $n = p_1^{a_1} p_2^{a_2} \cdots p_r^{a_r}$, 
$X_n$ is the conjunction $X_{p_1^{a_1}} \wedge X_{p_2^{a_2}} \wedge \cdots \wedge X_{p_r^{a_r}}$ of graphs $X_{p_1^{a_1}} , X_{p_2^{a_2}}, \cdots , X_{p_r^{a_r}}$, where 
a conjunction of graphs is defined as follows:

\begin{defn}\label{conjdef}
Let the graph $G_1$ have vertex set $V(G_1)$ and edge set $E(G_1)$, and graph $G_2$ have vertex set $V(G_2)$ and 
edge set $E(G_2)$. Then the \emph{conjunction} $G_1 \wedge G_2$ has vertex set 
$V(G_1 \wedge G_2) = V(G_1) \times V(G_2)$, and $(v_1, v_2)$ is adjacent to $(u_1, u_2)$ if $v_1u_1 \in E(G_1)$, and 
$v_2u_2 \in E(G_2)$.
\end{defn} 

Interestingly, our results can be extended to any conjunction $G_1 \wedge G_2 \wedge \cdots \wedge G_r$, where each 
$G_i$ is  complete $k_i$-partite.  Let $S = \{k_1, k_2, \dots, k_r\}$ be a multi-set of $r$ integers.  Let 
$\mathcal{G}^S = \{ G | G = G_1 \wedge G_2 \wedge \cdots \wedge G_r\}$, where $G_i$ is a complete $k_i$-partite graph.  
Denote the length of the longest induced cycle in $G \in \mathcal{G}^S$ by $\mathcal{M}(S)$, and
define $\mu(r) = \max_{S}{\mathcal{M}(S)}$ to be the length of the longest induced cycle in all graphs in $\mathcal{G}^S$, where $S$ contains 
$r$ integers.

\begin{thm}\label{conjthm}
For $r >1$, we have that $\mu(r) = m(r)$.
\end{thm}

To prove Theorem~\ref{conjthm}, we will create for conjunctions of $k_i$-partite graphs a representation similar to residue 
representation.  Then, using this representation, we will show how cycles in $G \in \mathcal{G}^S$ and $X_n$ are related.

\begin{defn}\label{partrep}
Let $S = \{k_1, k_2, \dots, k_r\}$, and let $G \in \mathcal{G}^S, G = G_1 \wedge G_2 \wedge \cdots \wedge G_r$.  Label the 
partitions in $G_i$ by $\{ 0, 1, 2, \dots, k_i-1\}$.  Let $v = (v_1, v_2, \dots, v_r) \in V(G)$, where $v_i$ 
belongs to 
partition $\alpha_i$ in $G_i$.  Then the \emph{partition representation} of $v$ is $\alpha_1 \alpha_2 \cdots \alpha_r$, 
and we say $v \simeq \alpha_1 \alpha_2 \cdots \alpha_r$.
\end{defn}

We can define the partition representation of a subgraph of $G \in \mathcal{G}^S$ as we defined the residue representation of a 
subgraph of $X_n$.  Namely, an induced subgraph on $\{x_1, x_2, \dots, x_l\}$ is written as an array of 
partition representations of the vertices $x_i$.  Note that an induced subgraph in $G$ is a cycle precisely when its 
partition representation satisfies the conditions needed for the residue representation of an induced cycle in 
$X_n$ -- no two non-consecutive rows can have similarities, and two non-consecutive rows must have at least one similarity.

\begin{proof}
[Proof of Theorem~\ref{conjthm}]

\item{(1)} First we show that $m(r) \geq \mu(r)$.  Suppose $S = \{k_1, k_2, \dots, k_r\}$, and $G \in \mathcal{G}^S$ contains an induced 
cycle $C$ of length $\mu(r)$, whose partition representation is 
\[
\begin{array}{cccc}
\alpha_{11}& \alpha_{12} &\cdots& \alpha_{1r}\\
\alpha_{21} &\alpha_{22}& \cdots &\alpha_{2r}\\
\vdots & \vdots & & \vdots\\
\alpha_{\mu(r)1} &\alpha_{\mu(r)2}& \cdots &\alpha_{\mu(r)r}\\
\end{array}.
\]
Note that, applying Proposition~\ref{nosame} 
to partition representations, no two rows above are identical if $\mu(r) > 4$.  So, if $\mu(r) > 4$, let 
$n = p_1 p_2 \cdots p_r$, where $p_i \geq \max{\{\alpha_{1i}, \alpha_{2i}, \dots, \alpha_{\mu(r)i}\}}$, 
and $p_i$ are prime.  Then the partition representation of $C$ 
above is in fact also the residue representation of an induced cycle in $X_n$, and so $X_n$ contains a cycle of 
length $\mu(r)$, as desired.  If $\mu(r) \leq 4$, we know that $m(r) \geq \mu(r)$, since $m(1) = 4$, and $m(r)$ increases with 
$r$ by Proposition~\ref{incr}.

\item{(2)} Now we show that $\mu(r) \geq m(r)$.  Let $X_n$, where  $n = p_1 p_2 \cdots p_r$, contain an induced cycle of length $m(r)$. 
Then $X_n \in G^{\{p_1, \dots, p_r\}}$, so $\mu(r) \geq m(r)$, as desired.
\end{proof}

Since our original problem concerns the Cayley graph $X_n$, we are also interested in Cayley graphs to which our results 
generalize.  In particular, we are interested in those graphs $G = Cay(A; A^*)$, where $A$ is a ring, $A^*$ 
is the group of units in $A$, and the graph $G$ is defined as follows:

(1) The vertex set $V(G)$ of $G$ is the set of elements in $A$.

(2) If $x, y \in V(G)$ then $\{x, y\} \in E(G)$, the edge set of $G$, if and only if $x-y \in A^*$.

We know that we can extend our results to any graph $G$ that is a conjunction of complete $k_i$-partite graphs 
for some $k_i$.  Note that, surprisingly, $k_i$ need not be finite, and, in fact, our Cayley graph need not 
contain a finite number of vertices.  For this, we rely on a partition using the Chinese Remainder Theorem.  One can 
refer to an algebra text such as \cite{Lang}, pp. 92-97 for the basic facts about rings and ideals needed to prove when such a 
partition gives us the desired graph structure.

\begin{defn}\label{locring}
A \emph{local ring} is a ring that contains only one maximal ideal.
\end{defn}

With this definition, we can show that a unitary Cayley graph on a product of local rings is a conjunction of complete 
$k_i$-partite graphs.

\begin{thm}\label{locrings}
Let $A_1, A_2, \dots A_r$ be local rings, and let $\mathfrak{m}_i$ be the one maximal ideal in $A_i$.  If 
$A = A_1 \times A_2 \times \cdots \times A_r$, then the Cayley graph $Cay(A; A^*)$ is a conjunction of complete 
$k_i$-partite graphs, for some nonzero $k_i$.
\end{thm}

\begin{proof}
[Proof of Theorem~\ref{locrings}]

We partition each ring $A_i$ into the $k_i$ residue classes modulo $\mathfrak{m}_i$, where $k_i = \#(A_i/\mathfrak{m}_i)$.  
Then, with this partition, we can show that the Cayley graph $Cay(A_i; A_i ^*)$ is a complete $k_i$-partite graph.  
Namely, $x, y \in A_i$ belong to the same residue class modulo $\mathfrak{m}_i$ if and only if $x-y \in \mathfrak{m}_i$ 
and is thus not a unit.  If $x, y \in A_i$ are in different residue classes modulo $\mathfrak{m}_i$, then 
$x-y \not \in \mathfrak{m}_i$.  Since $\mathfrak{m}_i$ is the only maximal ideal in $A_i$, and every non-unit element 
is contained in a maximal ideal of $A_i$, we have that $x$ and $y$ belong to different parts if and only if 
$x-y \in A_i^*$.  So two vertices in this graph belong to different parts if and only if they are 
adjacent.  So, indeed, $Cay(A_i; A_i ^*)$ is complete $k_i$-partite, where $k_i = \#(A_i/\mathfrak{m}_i)$.

Now we can show that $Cay(A; A^*)$ is a conjunction of complete $k_i$-partite graphs.  
We can assign every element of $A$ to some residue class modulo 
$\mathfrak{m}_iA$, for all $i$.  We also know that if $x, y \in A$, then $x \equiv y \pmod{\mathfrak{m}_iA}$ 
if and only if $x-y \in \mathfrak{m}_iA$.  This means that two vertices of $Cay(A; A^*)$ are not adjacent if and only if 
they belong to the same residue class modulo $\mathfrak{m}_iA$ for some $i$.  However, we can show that 
$x \not\equiv y \pmod{\mathfrak{m}_iA}$ for all $i$ if and only if $x-y \in A^*$.  Note that an element 
$z = (z_1, z_2, \dots , z_r) \in A$ is a unit in $A$ if and only if, for all $i$, $z_i \in A_i$ is a unit in $A_i$.  
So, since $z_i \in A_i$ is a unit if and only if $z_i \in \mathfrak{m}_iA$, we have that 
$x \not\equiv y \pmod{\mathfrak{m}_iA}$ for all $i$ if and only if $x-y$ is a unit in $A$.  So, indeed, $x, y \in A$ 
are adjacent if and only if they belong to different residue classes modulo $\mathfrak{m}_iA$ for all $1 \leq i \leq r$, 
and so $Cay(A; A^*)$ is a conjunction of complete $k_i$-partite graphs, where $k_i = \#(A/\mathfrak{m}_iA)$, as desired.
\end{proof}

Theorem~\ref{locrings} lets us extend our results to various unitary Cayley graphs.  In particular, it allows us to 
generalize to unitary Cayley graphs on Dedekind rings.

\begin{defn}\label{DD}
A \emph{Dedekind domain} (\cite{Marcus}) is an integral domain $R$ such that

(1) Every ideal in $R$ is finitely generated;

(2) Every nonzero prime ideal is a maximal ideal;

(3) $R$ is integrally closed in its field of fractions $$K = \{\alpha / \beta: \alpha, \beta \in R, \beta \not = 0\}.$$
\end{defn}

A Dedekind ring is simply a quotient of a Dedekind domain.


If $R$ is a Dedekind domain, and $\mathfrak{m}_i$ is a maximal ideal of $R$, then $R/{\mathfrak{m}_i}$ is a field and 
thus contains only one maximal ideal, $(0)$, and $R/{\mathfrak{m}_i^{a_i}}$ contains only the maximal ideal 
$\mathfrak{m}_i$, so $R/{\mathfrak{m}_i^{a_i}}$ is a local ring.  This is essential for the following corollary.

\begin{cor}\label{DR}
Let $R$ be a Dedekind domain, and let $I = \mathfrak{m}_1^{a_1}\mathfrak{m}_2^{a_2} \cdots \mathfrak{m}_r^{a_r}$ be a 
nonzero, non-unit ideal in $R$, where $\mathfrak{m}_i$ are maximal ideals of $R$.  Then the Cayley graph $Cay(A; A^*)$ 
is a conjunction of complete $k_i$-partite graphs, for $k_i = \#(R/{\mathfrak{m}_i})$.
\end{cor}

\begin{proof}

Since  $\mathfrak{m}_i$ are the distinct maximal ideals, $\mathfrak{m}_i^{a_i} + \mathfrak{m}_j^{a_j} = R$ for all 
$1 \leq i < j \leq r$.  Then the Chinese Remainder Theorem implies that 

$$A = R/{\mathfrak{m}_1^{a_1} \mathfrak{m}_2^{a_2} \cdots \mathfrak{m}_r^{a_r}} = R/{\mathfrak{m}_1^{a_1}} \times R/{\mathfrak{m}_2^{a_2}} \times \cdots \times R/{\mathfrak{m}_r^{a_r}}$$

We have noted above that $R/{\mathfrak{m}_1^{a_1}}$ is local, and thus we have that $A$ is a product of local rings.  By 
Theorem~\ref{locrings}, we have that the Cayley graph $Cay(A; A^*)$ is a conjunction of complete $k_i$-partite graphs, for 
$k_i = \#(R/{\mathfrak{m}_i})$.
\end{proof}

So, indeed, our theorems concerning $m(r)$ generalize to the maximum lengh of a cycle in unitary Cayley graphs on a 
Dedekind domain quotiented by an ideal with $r$ distinct maximal factors.  Dedekind domains are exactly those integral 
domains in which every ideal has a unique factorization into prime ideals, and thus are the rings of number theoretical 
interest.  Some nice examples of the Dedekind rings that we have generalized to above are the Gaussian integers 
modulo $a + bi$, denoted by $\mathbb Z[i]/{(a + bi)}$; any quotient of the ring of algebraic integers in the $p$th 
cyclotomic field $\mathbb Z[\zeta_p]$, where $\zeta_p$ is a $p$th root of unity; and any quotient of 
$\mathbb C[x, y]/{(y^2 - x^3 +x)}$, the ring of regular functions on the elliptic curve $y^2 = x^3 -x$.  Note that we 
also have generalized to unitary Cayley graphs on quotients of principal rings.

\section{Open Questions}
With the help of a computer program, written by Geir Helleloid, that performed an exhaustive search of 
arrays representing induced cycles, we have also been able to form conjectures about the lengths of the longest 
induced cycles in $X_n$.

The implementation of residue representation seems to promise more important results, both about the graph $X_n$, and 
more generally about conjunctions of complete $k_i$-partite graphs.  We know that the number of residues one can use 
to obtain a cycle of a given length $l$ is certainly bounded.  For example, the size of the residue set for a $6$-cycle 
for $r = 2$ cannot be greater than $3$.  In fact, if we can bound the size of the residue set needed to construct a cycle 
of length $m(r)$ to $a$, then $m(r) \leq a^r - (a-1)^r +2$, since the total number of possible vertices using $a$ 
residues is $a^r$, but every vertex among these vertices is adjacent to $(a-1)^r$ vertices, while in a cycle we want 
every vertex to be adjacent to exactly two other vertices.

Furthermore, we can continue reducing this bound of $a^r - (a-1)^r +2$, since among $a^r - (a-1)^r +2$ 
vertices there are still too many adjacencies for an induced cycle, and, in particular, there are too many vertices 
whose residue representation either contains no $1$'s or no $0$'s (see discussion in Section~\ref{ResRep}).

The computer program that we used to help predict the lower bound also seems to suggest that, not only may it be possible 
to modify any induced cycle to one of the same length whose residue set has size $3$, but that in fact we have 

\begin{conj}\label{tight}
$m(r) = 2^r + 2$.
\end{conj}

Actually, the computer program terminates for $r = 2$ and for $r = 3$, giving us that $m(2) = 6$ and 
$m(3) = 10$.  The program also gives us that the longest induced cycle one can construct using a residue set of only 
three residues for $r = 4$ has length $18$.  So, a question one may ask in verifying Conjecture~\ref{tight} is whether 
the longest cycle that uses only $3$ residues has length $2^r +2$.

There are several ways to approach these questions.  The most intuitive is to modify cycles of given lengths to 
cycles of the same length that use fewer residues.  However, we have not been able to find a general way of doing this 
for arbitrary cycles.  Another possibility is to show that any cycle can be modified to one of the same length 
that contains a column of only $2$ residues.  If so, we may ask whether we can reduce such a cycle in $X_N$, where $N$ has 
$r$ prime divisors, to a cycle in $X_n$, where $n$ has $r-1$ prime divisors, by deleting this column and a few rows to 
make the cycle induced.  Although we have yet to prove this, it seems that this method gives us a way of reducing an 
induced $k$-cycle in $X_N$ to an induced cycle of length approximately $k/2$ in $X_n$.  Since we know that, say, 
$m(3) = 10$, this could show that $m(r) \lesssim 10(2^{r-3})$ by induction.


Finally, one may also ask whether the use of residue representation can extend to graphs that are not conjunctions of 
complete $k_i$-partite graphs, and, if so, what conditions are necessary for our results to hold.

\section{Acknowledgements}
This research was done at the University of Minnesota Duluth. I would
like to give special thanks to Joseph A. Gallian for his encouragement and
support, and I would like to thank Philip Matchett and Melanie Wood, and Justin Sinz for their
insightful comments. I would also like to thank Geir Helleloid for helpful conversations. 
Funding was provided by grants form the NSF (DMS-0137611) and NSA (H-98230-04-1-0050).

\section{Appendix: Proof of conjecture 5.1}

\subsection{Alon's Theorem}

In \cite{Alon}, Noga Alon proved the following theorem:

\begin{thm}
Let $X_1, \ldots, X_n$ be disjoint sets, $r_1, \ldots, r_n, s_1, \ldots, s_n$ positive integers. For $1 \leq j \leq h$, let $A_j, B_j$ be subsets of $X := \coprod X_i$ such that

(1) $|A_j \cap X_i| \leq r_i, |B_j \cap X_i| \leq s_i,$ all $1 \leq i \leq n, 1 \leq j \leq h$.

(2) $A_i \cap B_i= \varnothing$

(3) $A_i \cap B_j \neq \varnothing , 1 \leq i < j \leq h$

Then $h \leq \prod \binom{r_i+s_i}{r_i}.$

\end{thm}

Although we will need only a very special case of this result, the proof of the general case is sufficiently short, important, and enlightening that it is reproduced here.

\begin{proof}
We may assume (by adjoining elements to each $A_i$ and $B_i$ subject only to the requirement that $A_i$ and $B_i$ remain disjoint) that $|A_j \cap X_i|=r_i$ and $|B_j \cap X_i|=s_i$. Let $V_i:=\mathbb{R}^{ r_i+s_i }$. For each $i$, choose a collection of vectors $\{ z_{i,t}|t \in ( \cup_j A_j \cup \cup_j B_j) \cap X_i \} $ in $V_i$ (that is, vectors indexed by elements in the above set) in general position. Define vector spaces $$V=\bigwedge_i( V_i^{\wedge r_i})$$ $$ \bar V =\bigwedge_i( V_i^{\wedge s_i});$$ both are subspaces of $\bigwedge (V_1 \oplus \ldots \oplus V_n)$ of dimension $\prod \binom{r_i+s_i}{r_i}$. Define elements $$y_j:=\bigwedge_i (\bigwedge_{t \in A_j \cap X_i} z_{i,t}) \in V$$
$$\bar y_j :=\bigwedge_i (\bigwedge_{t \in B_j \cap X_i} z_{i,t}) \in \bar V.$$  We have $y_i \wedge \bar y_i \neq 0, y_i \wedge \bar y_j=0, 1 \leq i < j \leq h$, just by properties of the wedge product (combined, of course, with the hypotheses on the intersections of the sets and the fact that the $z_{i,t}$ were chosen to be in general position); the wedges are scalars because $|A_j \cap X_i|=r_i$ and $|B_j \cap X_i|=s_i$. Since the matrix $(y_i \wedge \bar y_j)$ is invertible over $\mathbb{R}$, both the $\{ y_i \}$ and the $\{ \bar y_i \}$ are linearly independent subsets of $V$ and $\bar V$, respectively. Hence $h \leq \prod \binom{r_i+s_i}{r_i}.$
\end{proof}

This theorem naturally generalizes Bollob\'as' theorem \cite{Bol}, and its proof naturally generalizes Lov\'asz's exterior algebra proof \cite{Lo}. Note that when $r_i=s_i=1$ we obtain that $h \leq 2^n$.

\subsection{Application to Conjecture 5.1}
It was conjectured above that the longest induced cycle in the unitary Cayley graph $X_n = Cay(\mathbb Z_n , U_n)$ has length $2^r+2$, where $r$ is the number of prime divisors of $n$ (provided that either $r$ is greater than 1 or $r=1$ and the power of the prime is at least 2). Alon's theorem will be used to prove that $2^r+2$ is an upper bound, which when combined with Theorem 3.1 proves the conjecture.

\begin{proof}
[Proof of Conjecture 5.1]
The case of $r=1$ is treated separately in the above paper. Furthermore, the proof for the general case given here doesn't apply verbatim to the case of $r=1$ because of the failure of injection into the residue representation.

So fix $r \geq 2$; then by Theorem 2.10 it suffices to consider the case when $n=p_1 \ldots p_r$, a product of distinct primes. Let $v_1, \ldots, v_k$ be vertices of an induced cycle in $Cay(\mathbb Z_{p_1 \ldots p_r} , U_{p_1 \ldots p_r})$. Then the $v_i$ are represented by vectors in $\mathbb{N}^r$ via the residue representation, and distinct vertices correspond to distinct vectors by the Chinese Remainder Theorem. Each $v_i$ therefore gives a \emph {subset} $W_i$ of $\coprod_{1 \leq i \leq r} \mathbb{N}_i$ ($r$ copies of $\mathbb{N}$ each labeled by some $i$) in the natural way, and distinct vertices correspond to distinct subsets. The condition that the ${v_i}$ form an induced cycle implies that $$ W_i \cap W_{i+1} = W_i \cap W_{i-1}= \varnothing$$ and $$W_i \cap W_j \neq \varnothing, j \not \equiv i \pm 1 \pmod{k}$$  (see the discussion about similarity in and following Definition 2.2).

Let $$A_i := W_i, 1 \leq i \leq k-2,$$ $$B_i := W_{i+1}, 1 \leq i \leq k-2.$$ Then $|A_j \cap \mathbb{N}_i|=|B_j \cap \mathbb{N}_i|=1$, all $i,j$, $$A_i \cap B_i =W_i \cap W_{i+1} =\varnothing$$ and $$A_i \cap B_j =W_i \cap W_{j+1} \neq \varnothing, 1 \leq i < j \leq k-2,$$ because $(j+1) \equiv i+1 \Rightarrow j=i$ and $(j+1) \equiv i-1 \Rightarrow j+2 \equiv i$, but $1 < j+2 \leq k$ and $j+2 > j > i$. Therefore, by Alon's theorem, we have that $k-2 \leq 2^r$, i.e., $k \leq 2^r+2.$
\end{proof} 

\subsection{Acknolwedgements}
I would like to thank Josh Greene for pointing out some moments of substandard and misleading writing.

\end{document}